\documentclass[11pt]{amsart}

\usepackage{amssymb}
\usepackage{amsmath}
\usepackage{amsfonts}
\usepackage{amsthm}
\usepackage{amssymb}
\usepackage[mathscr]{eucal}
\usepackage[cp1250]{inputenc}

\newtheorem{theorem}{Theorem}[section]
\newtheorem*{theoremm}{Main Result}
\newtheorem{lemma}[theorem]{Lemma}

\newtheorem{remark}[theorem]{Remark}
\newtheorem{definition}[theorem]{Definition}

\begin{document}
\title{On submaximal curves in fake projective planes}

\author{Piotr Pokora and Halszka Tutaj-Gasi\'nska}
\address{Department of Mathematics, Pedagogical University of Cracow,
	Podchor\c a\.zych 2,
	PL-30-084 Krak\'ow, Poland.}
\email{piotrpkr@gmail.com, piotr.pokora@up.krakow.pl}
\address{Jagiellonian Univeristy, Faculty of Mathematics and Computer Science, {\L}ojasiewicza 6, PL-30-348 Cracow Poland}
\email{Halszka.Tutaj@im.uj.edu.pl}

\begin{abstract}
The main purpose of the paper is to exclude the existence of certain submaximal curves in fake projective planes. This will lead to lower bounds on multipoint Seshadri constants of `fake' $\mathcal{O}(1)$ on fake projective planes.
\textbf{AMS2010}: 14J29, 14C20
\end{abstract}

\maketitle

\section{Introduction}
In the paper, we study the existence of certain submaximal curves in the context of the Seshadri constants for ample line bundles on fake projective planes. Let us recall that by a fake projective plane we understand a smooth complex projective surface of general type having the same Betti numbers as the complex projective plane. The existence of such fake projective planes was proved by Mumford in \cite{mum}, and now we know that there are exactly $50$ pairs of fake projective planes, see \cite{Cart}. We would like to focus on the case of multipoint Seshadri constants for fake projective planes. In the case of the single point Seshadri constants, L. Di Cerbo \cite{Luca} proved that these constants coincide with the single point Seshadri constants of the complex projective plane. The key advantage of his results is that it provides (probably) the first sharp result on single point Seshadri constants for surfaces of general type. Here we want to follow this path and look at the multipoint Seshadri constants (that are defined precisely at the beginning of Section \ref{prel}) for ample line bundles on fake projective planes. 

Let $X$ be a fake projective plane, i.e., a surface with $b_{1}(X) = 0, b_{2}(X) = 1$ and $X$ is not isomorphic to $\mathbb{P}^{2} = \mathbb{P}^{2}_{\mathbb{C}}$. Then its canonical bundle is ample and $X$ is of general type. So a fake projective plane is nothing but a surface of general type with $p_{g} = 0$ and $c_{1}^{2}(X) = 3c_{2}(X) = 9$. Furthermore, its fundamental group $\pi_{1}(X)$ is infinite. Moreover, $X$ has Picard number $1$, and we can find an ample generator $L_{1}$, which turns out to be not uniquely determined, such that $L_{1}^2 = 1$. 
 
There are not so many results on multipoint Seshadri constants for surfaces of general type and  ample line bundles (in some places the assumption on the ampleness is weakened to big and nef line bundles). In \cite{Szemberg08}, T. Szemberg showed that if $X$ is a smooth complex projective surface of Picard number $1$ and $L$ is the ample generator of the N\'eron-Severi group of $X$, then
$$\varepsilon(X,L;r) \geq \bigg\lfloor \sqrt{\frac{L^{2}}{r}} \bigg\rfloor.$$

In the same context, one can show that if $X$ is a smooth complex projective surface of Picard number $1$, $L$ is the ample generator of the N\'eron-Severi group, and the number of points $r=k^{2}$ for some $k\in \mathbb{Z}_{\geq 1}$, then
$$\varepsilon(X,L;r) = \frac{\sqrt{L^{2}}}{k},$$
so in the present paper \emph{we are going to focus on the non-trivial case when $r\neq k^{2}$ for some $k$}.

As it was observed by J. Ro\'e \cite{Roe}, there is an interesting link between multi-- and single point Seshadri constants, so let us phrase \cite[Theorem 3]{Roe} in the setting of our paper and adequate for the scope of our work. If $X$ is a fake projective plane, $p\in X$ any point, and $r\geq 1$, then one has
$$\varepsilon(X,L_{1}; r) \geq \varepsilon(X, L_{1}; p) \cdot \varepsilon(\mathbb{P}^{2},  \mathcal{O}_{\mathbb{P}^{2}}(1);r).$$
By a result due to L. Di Cerbo \cite{Luca} we know that $\varepsilon(X,L_{1};p)=1$, so we obtain the following inequality
$$\varepsilon(X,L_{1}; r) \geq \varepsilon(\mathbb{P}^{2}, \mathcal{O}_{\mathbb{P}^{2}}(1);r).$$
The main result of our paper, which is Theorem \ref{main}, gives us bounds for Seshadri constants on fake projective planes. It tells us, for instance, that for $r \in \{2,3,5,6,7\}$ we have strict inequality above. This stands to the opposite to the case of single point Seshadri constants where $\varepsilon(X,kL_{1};P) = \varepsilon(\mathbb{P}^{2}, \mathcal{O}_{\mathbb{P}^{2}}(k); P')$. Our strategy to show the strict inequality is based on the fact that we are able to exclude the existence of certain \emph{submaximal} curves, i.e., those irreducible and reduced curves in fake projective planes breaking the so-called Nagata bound. We explain in detail these notions in the next section.

 The main reason for our studies on multipoint Seshadri constants for fake projective planes is the Nagata-Biran-Szemberg Conjecture, see \cite{Primer}. It claims that the multipoint Seshadri constant of an ample line bundle $L$ at a very general set of points is maximal when $r\geq k^{2}_{0}L^{2}$, where $k_{0}$ is the smallest integer such that the linear system $|k_{0} L|$ contains a smooth non-rational curve. In the case of fake projective planes we know that there are no rational and elliptic curves on them,  which means that multipoint Seshadri constants in this case must be maximal unconditionally.
 
We conclude the Introduction with a somewhat less technical version of our main result, which in full detail is stated in Theorem \ref{main}.
 \begin{theoremm}
	Let $X$ be a fake projective plane and denote by $L_{1}$ an ample generator of the N\'eron-Severi group.	Then the multiple Seshadri constant of $L_{1}$ at $r$ (very) general points with $r\geq 2$ is bounded from below by  $$\varepsilon(X, L_1 ;r)\geq \frac{1}{\sqrt{r}+0.31}.$$
\end{theoremm}	
\section{Preliminaries}
\label{prel}
Let us start with the definitions of single and multipoint Seshadri constants.
\begin{definition} Let $X$ be a complex smooth projective surface and $L$ an ample line bundle. For a point $x \in X$ the Seshadri constant of $L$ at $x$ is defined as
	$$\varepsilon(X,L;x) = \displaystyle{{\rm inf}_{C \ni x}} \frac{C.L}{{\rm mult}_{x}C},$$
	where the infimum is taken over all irreducible and reduced curves $C$ passing through $x \in X$.
\end{definition}
\begin{definition} Let $X$ be a complex smooth projective surface and $L$ an ample line bundle at $r$ points $x_{1}, ..., x_{r}$, then the multipoint Seshadri constant of $L$ is defined as
	$$\varepsilon(X, L; x_{1}, ..., x_{r})=\inf\bigg \{\frac{C.L}{\sum_{i=1}^r \text{mult}_{x_i}C} \bigg\},$$
	where the infimum is taken over all irreducible and reduced curves such that $C \cap \{x_{1}, ..., x_{r}\} \neq \emptyset$. If the points $x_{1}, ..., x_{r}$ are very general, then we use the abbreviation $\varepsilon(X, L;r)$.
	\end{definition}
Let $x_1,\dots,x_r$ be (very) general points in a surface $X$ (which is always smooth, projective and complex). Let us recall that by K. Oguiso's result \cite{Oguiso} the maximum of the multipoint Seshadri constants is achieved for a generic set of $r$ points. That is why we  use the abbreviation $\varepsilon(X,L;r)$ for the $r$-point Seshadri constant with $r$ points in very general position. Let $C$ be a curve on $X$ and  denote $m_i=\text{mult}_{x_i}C$. We recall some fundamental results about curves by G. Xu \cite[Lemma 1]{xu} and M. Roth \cite{Roth}.

\begin{lemma}[Xu] \label{xu}
Let $X$ be a smooth complex projective surface and $C_{U}=\{(C_{u},x_{u}) \, | \,u\in U\}$ be a non-trivial family of pointed irreducible and reduced curves in X such that ${\rm mult}_{x_{u}} C_{u} \geq m$ for some integer $m\geq 2$, then
$$C^{2}_{u} \geq m(m-1) + 1.$$
\end{lemma}
\begin{remark}\label{xu-improved}{\rm
In our consideration, we will use the result by A. L. Knutsen, T. Szemberg, W. Syzdek and also by F. Bastianelli,  see \cite{KSS} and \cite{Bas}, extending the result of Xu, saying that
	for a smooth projective surface $X$ and  a smooth variety $U$, and a nontrivial family $\{(C_u,x_u)\}_{u\in U}$, where $x_u$ is a very general point of $X$ and  $C_u$ is a curve satisfying  the condition $\text{mult}_{x_u}C_u \geq m$ for every  $u\in U$ and for some integer $m\geq 2$, we have for a general curve $C$ of this family 
	$$C^2 \geq m(m-1)+\text{gon}(\widetilde{C}),$$
	where $\widetilde{C}$ denotes the normalization of $C$.
 As mentioned above, on fake projective planes there are no rational or elliptic curves, so $\text{gon}(\widetilde{C})\geq 2$, and we will use the inequality
   \begin{equation}\label{xu-imp}
   C^2 \geq m(m-1)+2.
    \end{equation}}
\end{remark}

In order to state next useful results we need the following definition.
\begin{definition}
	Let $X$ be a smooth projective surface and $L$ an ample line bundle. The Seshadri constant is said to be submaximal if the inequality holds
	$$\varepsilon(X,L; r) < \sqrt{\frac{L^{2}}{r}}.$$
	In that case a curve $C \subset X$ which computes the ratio $\frac{C.L}{\sum_{i}{\rm mult}_{x_{i}}C}$ is called \emph{submaximal}.

\end{definition}
We have the following important results for submaximal curves which are proven by M. Roth \cite{Roth} (manuscript on request).

\begin{theorem}\label{roth}
Let $X$ be a smooth complex projective surface and $C \subset X$ a submaximal curve with respect to a fixed $L$ and very general points $\{x_{1}, ..., x_{r}\}$ with $r\geq 1$. Then we have the following properties:
\begin{enumerate}
	\item[a)] At least $(r - 1)$ of the $m_{i}$'s are equal.
	\item[b)] If $m_i$ and $m_{j}$ are distinct multiplicities, then $$-D^{2} \leq (m_{i} - m_{j})^{2} < -\frac{r}{r-1}D^{2},$$ where $D$ is the proper transform of $C$ under the blowing up of $X$ at $\{x_{1}, ..., x_{r}\}$.
	\item[c)] If one of the multiplicities is equal to zero, then others non-zero multiplicities are equal to $1$ and $C^{2} = -1$.
	\item[d)] If none of the multiplicities is zero and $(r,m_{r}) \neq (1,1)$, then $\sum_{i=1}^{r} m_{i} = \lceil \sqrt{rC^{2}} \rceil$.
	\item[e)] In the case when all multiplicities are non-zero and additionally equal, say to $m \geq 1$, then
	$$m = \bigg\lceil \sqrt{\frac{C^{2}}{r}} \bigg\rceil, \quad m - \sqrt{\frac{C^{2}}{r}} \leq \frac{1}{r}, \quad \text{ and } \quad rm^{2} - C^{2} \leq m.$$
	\item[f)] If none of the  multiplicities are zero and one of the multiplicities is different from the others, then
	$$\sum_{i=1}^{r} m_{i} - \sqrt{rC^{2}} < \frac{1}{r}.$$
\end{enumerate}
\end{theorem}
\begin{remark}
It is worth pointing out that Theorem \ref{roth} a) is a nice generalization of Corollary 4.6 from T. Szemberg's Habilitation Thesis \cite{Szemberg} where he was considering polarized surfaces $(X,L)$ having Picard number $\rho(X) = 1$.
\end{remark}
\begin{remark}
It is worth emphasizing that Theorem \ref{roth} a) was also proved by B. Harbourne and J. Ro\'e in \cite[Corrolary 2.2.2]{HR}.
\end{remark}

Now we focus on fake projective planes.

\begin{remark}
For  a  fake  projective  plane  $X$,  we  have $Pic(X) \simeq H^{2}(X; \mathbb{Z})$. Moreover,  the  torsion  part  of $H^{2}(X;\mathbb{Z})$ is  equal  to $H_{1}(X;\mathbb{Z})$ which is never vanishing.
\end{remark}
\begin{remark}
For a fake projective  plane  $X$,  we  denote  by $L_{1}$ any ample generator of $Pic(X)$ modulo torsion. Let us emphasize that $L_{1}$ is not uniquely determined.

\end{remark}
Here we are going to find lower bounds on the multipoint Seshadri constant for fake projective planes and a specially chosen line bundle, namely $L_{1}$. Let us recall that $L_{1}$, being ample, has $L_{1}^{2} = 1$, but it is an open problem to determine whether in general one has $h^{0}(X,L_{1}) > 0$. Here we are not going to approach this question, but it seems to be very interesting to check whether $L_{1}$ might be effective due to important implications -- please consult \cite[Remark 3.9]{Luca1}.

\section{Main Result}

The bound on the Seshadri constants in $r$ general points on fake projective plane is given by the following theorem. 

\begin{theorem}\label{main}
	Let $X$ be a fake projective plane and denote by $L_{1}$ an ample generator of the N\'eron-Severi group and assume that $r$ is not a square. Then $$\varepsilon(X, L_1 ;r)\geq \frac{1}{\sqrt{r}+\delta(r)},$$
		where
		$$\delta(2)=0.031, \quad \delta(3)=0.018, \quad \delta(5)=0.014,$$
		$$ \delta(6)=0.022, \quad \delta(7)=0.011, \quad \delta(8)=0.012,$$
		 and 
		 	$$ \delta(r)=0.013 \text{ for } r\geq 10.$$
		 	If  $r=s^2$ for $s\in \mathbb{Z}_{>0}$, then 
		 	$$\varepsilon(X, L_1 ;r) = \frac{1}{s}.$$
	\end{theorem}

Before we present our proof of the Theorem, let us make a few remarks.

\begin{remark}\label{szembergbound}
{\rm
\begin{itemize}
	\item[1)] For $\mathbb{P}^{2}$ we have the following lower bound \cite[Theorem 2.3]{SzSz}: 
	$$\varepsilon(\mathbb{P}^{2}, \mathcal{O}_{\mathbb{P}^{2}}(1); r)\geq \frac{\sqrt{49r+8}}{7r+1}$$
	provided that $r\geq 10$.
	Thus, combining \cite[Theorem 3]{Roe} with the above lower bound we get 
	$$\varepsilon({\rm FPP}, L_{1}; r) \geq \frac{\sqrt{49r+8}}{7r+1},$$
	provided that $r\geq 10$. 
	\item[2)]	It is easy to check that 
	$$\frac{1}{\sqrt{r}+0.013}>\frac{\sqrt{49r+8}}{7r+1}$$
	for $10\leq r\leq 22,$
	so our bound is better in this range.
		\item[3)] Analysing the proof of our theorem,  the reader may notice that changing the range of $r$ to improve 
		the bound, e.g. for $23 \leq r$  we  have
		$$\varepsilon({\rm FPP}, L_{1}; r) \geq \frac{1}{\sqrt{r}+0.010},$$
	  so having $r$ big enough (and a strong enough computer) it may be possible to
	  find really good bounds for $\varepsilon({\rm FPP}, L_{1}; r)$.
\end{itemize}
 }

\end{remark}

\begin{remark}
	Here we would like to compare multipoint Seshadri constants for fake projective planes and the projective plane.
	\begin{center}
		\begin{tabular}{c|c|c}
			$r$	& $\varepsilon(\mathbb{P}^{2}, \mathcal{O}_{\mathbb{P}^{2}}(1);r)$  & $\varepsilon({\rm FPP}, L_{1};r)$ \\ \hline \hline 
		     	& Exact value & Bound from Thm \ref{main} \\
			\hline 
			$2$	 &$ 1/2$ &  $> 0.69$      \\ 
			$3$	 & $1/2$  &  $\geq 0.5701$ \\ 
			$4$	 & $1/2$  &  $ 1/2$         \\ 
			$5$	 & $2/5$  &  $\geq 0.44 $   \\ 
			$6$	 & $2/5$  &  $\geq 0.4046$ \\ 
			$7$	 & $\frac38=0.375$  &  $\geq 0.3763$ \\ 
			$8$	 & $\frac{6}{17}\approx 0.3529$ &  $\geq  0.3391 $  \\
			$9$  & $1/3$  &   $1/3$  \\
			\hline
			    & Bound from Remark \ref{szembergbound} & Bound from Thm \ref{main} \\
			\hline 
			$10$ & $\geq 0.3143$       &    $\geq 0.3149$    \\
			$11$ &  $\geq 0.2998$      &    $\geq 0.3003$    \\
			$12$ &  $\geq 0.2872 $      &    $ \geq 0.2876$    \\
			$13$ &  $\geq0.2760 $      &   $ \geq 0.2763$     \\
			$14$ &  $\geq 0.2661 $      &    $ \geq 0.2663$    \\
			$15$ &  $\geq 0.2571 $      &       $ \geq 0.2573$ \\
			$16$ &  $1/4$     &      $1/4$  \\
		\end{tabular} 
	\end{center}
\end{remark}

Now we present the proof of Theorem \ref{main}. 

\begin{proof} 
	The idea of the proof goes as follows. Suppose we have a submaximal curve passing in fake projective plane $X$ through $r$ points with multiplicities $m$ (in $r-1$ points) and $M$ (in one point) -- it is also possible that $M=m$.  Using Theorem \ref{roth} we check that the inequality $\frac{k}{(r-1)m+M}\geq \frac{1}{\sqrt{r}+\delta(r)}$ is satisfied for  $k$ big enough. The remaining cases (small $k$) are excluded with the help of the Xu-type statement (\ref{xu-imp}).

Let us now go into details.
	Suppose there exists a submaximal curve $C\in |kL_1|$ passing through  general points $x_1\dots x_r$ with multiplicities $m$  (in $r-1$ points) and $M$ (in one point).
Using M.~Roth's result, i.e., Theorem \ref{roth}, we have  that every submaximal curve must satisfy
	$$\frac{C.L_{1}}{(r-1)m+M}\geq \frac{k}{k\sqrt{r}+\frac{1}{r}}$$
	provided that $M\neq m$, and if $m=M$ we have
	$$\frac{C.L_{1}}{rm}\geq \frac{k}{k\sqrt{r}+\frac{1}{2}}.$$
	In both cases we have
	$$\frac{C.L_{1}}{(r-1)m+M}\geq \frac{k}{k\sqrt{r}+\frac{1}{2}}.$$
	Next, we check when
	$$\frac{k}{k\sqrt{r}+\frac12}\geq  \frac{1}{\sqrt{r}+\delta(r)}.$$
	This may be checked of course only for the smallest considered $\delta$, so 
	we need to check when 
		$$\frac{k}{k\sqrt{r}+\frac{1}{2}}\geq  \frac{1}{\sqrt{r}+0.01}.$$
	The inequality is satisfied for $k\geq 50$, so this means that all submaximal curves violating the bound $\frac{1}{\sqrt{r}+\delta(r)}$, if they exist, must be in $|kL_1|$ with $k\leq 49$.

	Firstly, observe that a submaximal curve must have at least one multiplicity greater than $1$. Indeed, if all points have multiplicity one, the conditions imposed by them are independent, thus computing $h^0(kL_1)$
and imposing $r$ conditions we get that 
$$k\geq \frac{3+\sqrt{1+8r}}{2}$$
whereas the submaximality gives
$$\frac{k}{r}< \frac{1}{\sqrt{r}},$$
a contradiction.
	Thus, at least one multiplicity of a submaximal $C$ is greater or equal $2$.
	So, from Inequality (\ref{xu-imp}), every submaximal curve $C$ must satisfy one of the following 
	inequalities:
	$$f_1(k,r,m):= rm^2-m+2-k^2\leq 0 \text{ if only  $M=m$,}$$
or	
	$$f_2(k,r,m,M):=(r-1)m^2+M^2-M+2-k^2\leq 0 \text{ if only $1<M<m$,}$$
	 or
	$$f_3(k,r,m,M):= (r-1)m^2+M^2-m+2-k^2\leq 0 \text{	if only $1<m<M$,}$$
or 
	$$f_4(k,r,m):= (r-1)m^2+1-m+2-k^2\leq 0 \text{ if $ m_r=1, m_1=\dots m_{r-1}=m>1$,}$$
		or 
	$$f_5(k,r,m):= (r-1)+m^2-m+2-k^2\leq 0 \text{ if $ m_1>1, m_2=\dots m_{r}=1$.}$$

	We only have to check that if $r\geq 2$ and $k\leq 49$ there are no positive integer solutions (for a given $k$, $m, M$, and $r$) for the systems of inequalities given below.
	If $m=M$, then
	$$f_1\leq 0, \  \frac{k}{rm}< \frac{1}{\sqrt{r}+\delta(r)}.$$
	If $M<m$, then
	$$f_2\leq 0, \  \frac{k}{(r-1)m+M}< \frac{1}{\sqrt{r}+\delta(r)}.$$
	In the case when $m<M$:
	$$f_3\leq 0, \  \frac{k}{(r-1)m+M}< \frac{1}{\sqrt{r}+\delta(r)}.$$
	Finally, when one multiplicity is $1$ and the others are greater than $1$:
	$$f_4\leq 0, \  \frac{k}{(r-1)m+1}< \frac{1}{\sqrt{r}+\delta(r)}$$
	and when all but one multiplicities are $1$:
		$$f_5\leq 0, \  \frac{k}{(r-1)+m}< \frac{1}{\sqrt{r}+\delta(r)}$$
	This may be done by hand, or, preferably with help of a computer program, e.g. Wolfram Mathematica.
\end{proof}

\begin{remark}{\rm 
	In order to get a better bound on $\varepsilon(X, L_1;2)$, namely
	$$\varepsilon(X, L_1;2)>0.7$$
	we would have to exclude, for example, the existence of 
	$C\in |7L_1|$ having exactly two singular points of multiplicity $m=5$,
what seems to be not obvious.	}
\end{remark}

\section*{Acknowledgments}
We would like to thank L. Borisov, L. Di Cerbo, and M. Stover for discussions on fake projective planes, and to Joaquim Ro\'e for useful suggestions. Additionally, we would like to warmly thank Mike Roth for sharing with us his manuscript and very useful discussion on the subject of the present note. Finally, we would like to thank an anonymous referee for suggestions that allowed to improved the readability of the note.
The first author is partially supported by the National Science Center (Poland) Sonata Grant \textbf{Nr 2018/31/D/ST1/00177}.

\end{document}